\documentstyle[12pt,epsf]{article}

\newcommand{\complex}{{\bf C}}

\newcommand{\integer}{{\bf Z}}

\newcommand{\cP}{{\cal P}}

\newtheorem{theorem}{Theorem}[section]
\newtheorem{numtheorem}{Theorem}[section]
\newtheorem{example}[theorem]{Example}
\newtheorem{proposition}[theorem]{Proposition}
\newtheorem{numproposition}[numtheorem]{Proposition}

\newtheorem{lemma}{Lemma}
\newtheorem{mlemma}{Main Lemma}
\newtheorem{corollary}[theorem]{Corollary}

\newcommand{\prf}{{\it Proof:} }
\newcommand{\qed}{\hspace*{\fill}\hbox{$\Box$}}
\newtheorem{numexample}{Example}

\begin{document}

\title{Quantum multiplication of Schur Polynomials}
\author{{\sc Aaron Bertram}\thanks{Department of Mathematics, University of Utah, Salt Lake City, UT 84112, {\em email address:}\ bertram@math.utah.edu },\and {\sc Ionu\c t Ciocan-Fontanine}\thanks{Institut Mittag-Leffler, Aurav\" agen 17, S-182 62, Djursholm, Sweden, {\em email address:} ciocan@ml.kva.se}, \and {\sc and William Fulton}\thanks{Institut Mittag-Leffler, Aurav\" agen 17, S-182 62, Djursholm, Sweden, and University of Chicago, Chicago, IL 60637, {\em email address:}\ fulton@math.uchicago.edu } }
\date { May 25, 1997}
\maketitle

\section{ Introduction}

The Giambelli formula writes the class of a Schubert variety in a Grassmannian as a Schur polynomial, and rules of Pieri and Littlewood-Richardson give explicit formulas for multiplying these classes.
The quantum cohomology ring of a Grassmannian is known, together with formulas for the classes of Schubert varieties \cite{B}. Our aim here is to present and discuss formulas for their quantum multiplication.

Fix positive integers $k$ and $l$, and set $n=k+l$. Let $X:=Gr(l,n)$ be the Grassmannian of $l$-planes in $\complex ^n$. Let $\sigma_1,\dots ,\sigma_k$ be indeterminates, with $\sigma_i$ of degree $i$.
Define polynomials $Y_r=Y_r(\sigma_1,\dots ,\sigma_k)$ by the formulas
\begin{equation}
{Y_r={\rm det}(\sigma_{1+j-i})_{1\leq i,j\leq r}.}
\end{equation}
In the classical cohomology ring $H^*(X)$, if $\sigma_i$ maps to the $i^{\rm th}$ Chern class of the universal quotient bundle, then $(-1)^rY_r$ maps to the $r^{\rm th}$ Chern class of the universal subbundle. This leads to one 
of the standard presentations of the classical cohomology ring of the Grassmannian:
\begin{equation}
H^*(X)=\integer [\sigma_1,\dots ,\sigma_k]/(Y_{l+1},\dots ,Y_n).
\end{equation}
The (small) quantum cohomology ring $QH^*(X)$ is an algebra over $\integer [q]$,
with $q$ a variable of degree $n$. It has the presentation (\cite{ST}, \cite{B}, \cite{FP}) 
\begin{equation}
QH^*(X)=\integer[q,\sigma_1,\dots ,\sigma_k]/(Y_{l+1},\dots ,Y_{n-1},Y_n+(-1)^kq).
\end{equation}

We identify a partition $\lambda$ with its Young diagram;
we write $\lambda\subset l\times k$ for a partition
$\lambda =(\lambda_1\geq\dots\geq\lambda_l\geq 0)$ with $\lambda_1\leq k$, i.e., for a partition whose Young diagram fits inside an $l\times k$ rectangle.

\vskip 10pt
\epsfxsize 9cm
\centerline{\epsfbox{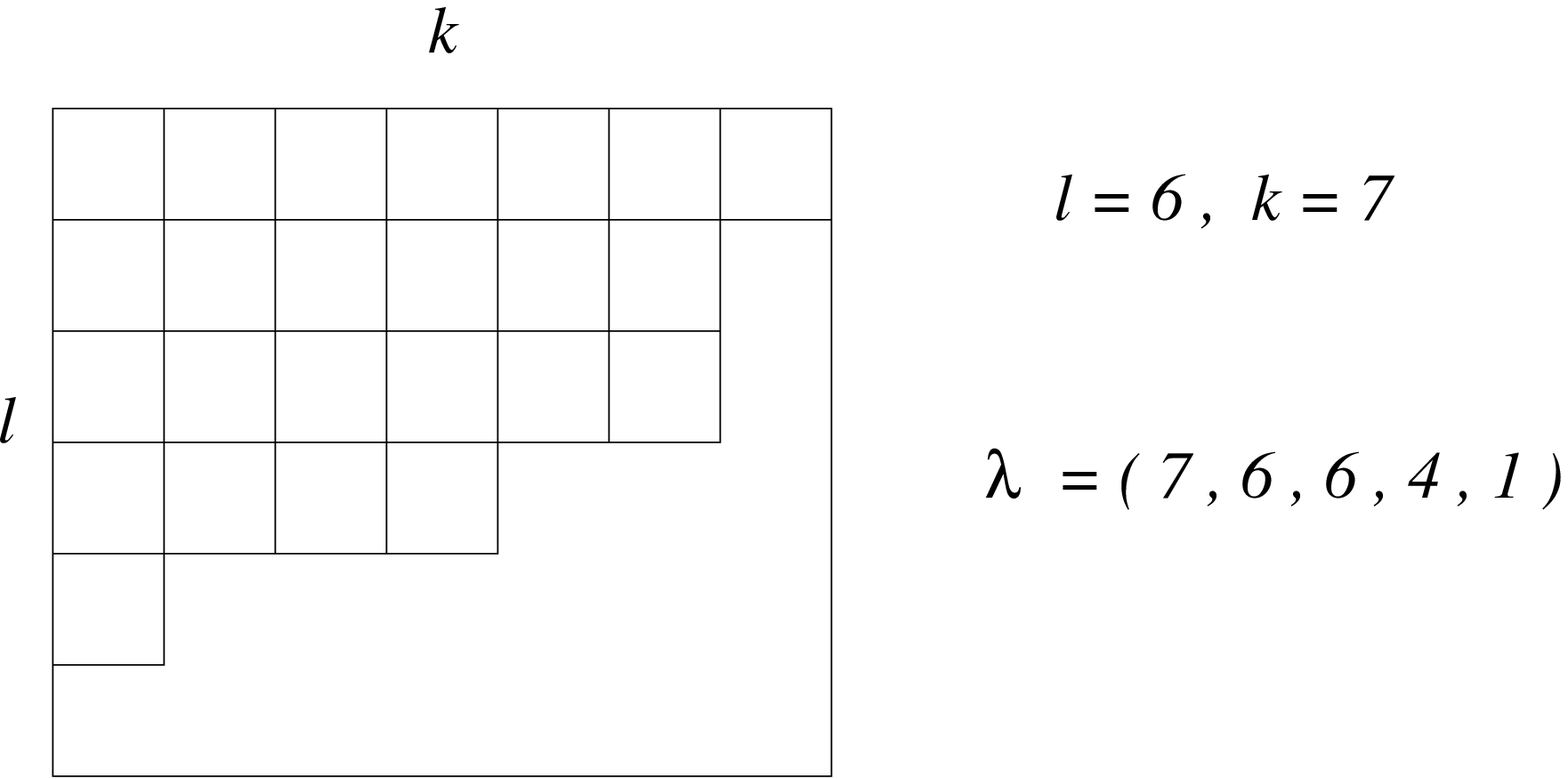}}
\vskip 20pt

There is a Schubert variety $\Omega_{\lambda}\subset X$ for each 
$\lambda\subset l\times k$; $\Omega_{\lambda}$ consists of $l$-dimensional
subspaces of $\complex^n$ such that
${\rm dim}( L\cap\complex^{k+i-\lambda_i}) \geq i$, for $1\leq i\leq l$.
The classical Giambelli formula says that the class of $\Omega_{\lambda}$ in
$H^*(X)$ is $\sigma_{\lambda}$, where
\begin{equation}
\sigma_{\lambda}={\rm det}(\sigma_{\lambda_i+j-i})_{1\leq i,j\leq l}.
\end{equation}
These classes $\sigma_{\lambda}$ give a $\integer$-basis for $H^*(X)$. Their product in this ring can therefore be written
\begin{equation}
\sigma_{\lambda}\cdot\sigma_{\mu}=\sum N_{\lambda\mu}^{\nu}\sigma_{\nu}.
\end{equation}
Here $\lambda,\mu$ and $\nu$ are partitions inside the $l\times k$ rectangle,
with sizes $ | \nu |  = | \lambda |  + | \mu | $.
That the coefficients $N_{\lambda\mu}^{\nu}$ are nonnegative follows from the fact that the Schubert varieties can be translated by the action of $GL_n(\complex)$ so that they meet properly. In fact, Grassmannians are one class of varieties for which explicit closed expressions are known for multiplying Schubert classes, which also show the coefficients to be nonnegative. The formula for multiplication is known as the {\em Littlewood-Richardson rule}, and states that $N_{\lambda\mu}^{\nu}$ is the
number of tableaux on the skew shape $\nu/\lambda$ of content $\mu$, whose word is a reverse lattice word. Here $\nu/\lambda$ is the complement of the Young diagram of $\lambda$ in the Young diagram of $\nu$. A {\em tableau} of {\em content} $\mu=(\mu_1,\dots ,\mu_l)$ is a numbering of the boxes of $\nu/\lambda$ with $\mu_1$ $1$'s,  $\mu_2$ $2$'s, $\dots$ , up to  $\mu_l$ $l$'s, which are weakly increasing across rows, and strictly increasing down columns (these are often called {\em semistandard tableaux}). The {\em word} of a tableau is the list of its entries, read from left to right in rows, from bottom to top. A word is {\em reverse lattice} if, from any point in it to the end, there are at least as many $1$'s as there are $2$'s, at least as many $2$'s as $3$'s, and so on. An example with $\lambda =(2,1)$, $\nu =(3,3,3,1)$, $\mu =(3,2,2)$ is shown in the figure below.

\vskip 20pt
\epsfxsize 55mm
\centerline{\epsfbox{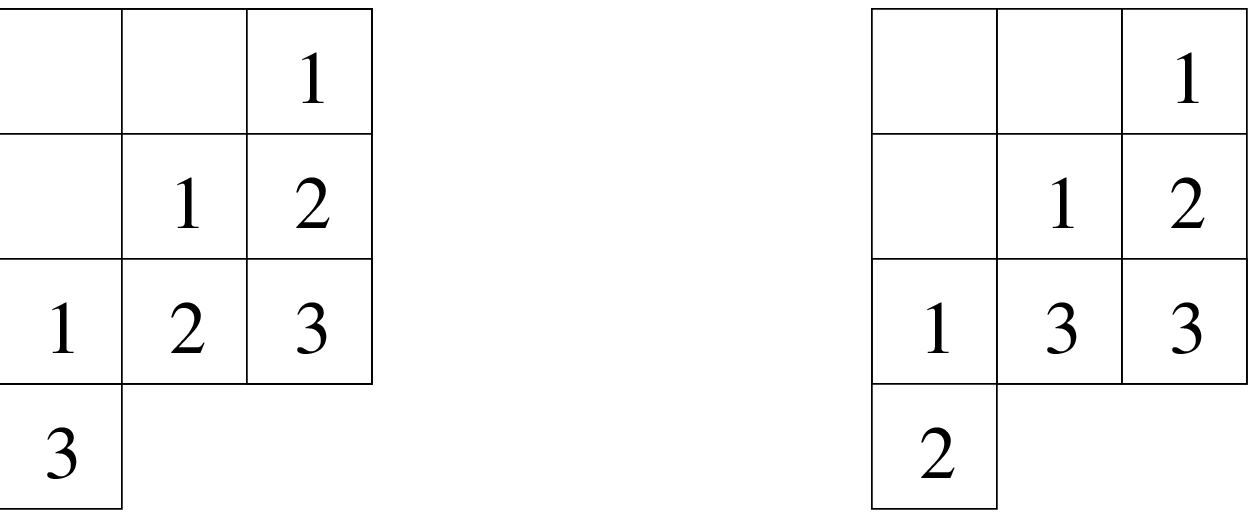}}
\vskip 20pt

The word $3123121$ of the first tableau is a reverse lattice word, while
the word $2133121$ of the second is not.

The rule follows from the corresponding identity for Schur functions
\begin{equation}
s_{\lambda}\cdot s_{\mu}=\sum N_{\lambda\mu}^{\nu}s_{\nu},
\end{equation}
where $s_{\lambda}={\rm det}(s_{\lambda_i+j-i})$, with $s_1,s_2,\dots $ indeterminates. In this case there is no need to restrict to partitions inside any rectangle; the map 
$$s_i\longmapsto \sigma_i\in H^*(X)$$ takes $s_{\lambda}$ to $0$ if $\lambda$ is not contained in the $l\times k$ rectangle (see \cite{M}, \cite{F2}).

The Pieri formula is a special case,
\begin{equation}
s_p\cdot s_{\lambda}=\sum s_{\nu},
\end{equation}
the sum over all $\nu$ obtained from $\lambda$ by adding $p$ boxes, with no two in the same column. Repeated, it gives a formula
\begin{equation}
s_{\mu _1}\cdot\dots\cdot s_{\mu _r}\cdot s_{\lambda}=\sum K_{\lambda\mu}^{\nu}s_{\nu}.
\end{equation}
Here $\mu _1,\dots ,\mu_r$ are nonnegative integers, the sum is over all $\nu$
with $ | \nu |  = | \lambda |  +\sum_{i=1}^r\mu_i$, and $K_{\lambda\mu}^{\nu}$ is the {\em Kostka number}: the number of tableaux on $\nu/\lambda$ with content $\mu$. The usual Kostka number $K_{\nu\mu}$ is $K_{\lambda\mu}^{\nu}$, where $\lambda=\emptyset$ and $\mu$ is a partition. It is a consequence of (8) that $K_{\lambda\mu}^{\nu}$ is independent of the order of terms in $(\mu _1,\dots ,\mu_r)$.

We turn now to the quantum analogues of these formulas. We identify $QH^*(X)$ with the ring presented in (3), and we regard $\sigma_1,\dots ,\sigma_k$ as elements of this ring. For any partition $\lambda =(\lambda_1,\dots ,\lambda_r)$, we let 
$$\sigma_{\lambda}={\rm det}(\sigma_{\lambda _i +j-i})_{1\leq i,j\leq r},$$ 
an element of $QH^*(X)$.
The classes $\sigma_{\lambda}$, for $\lambda\subset l\times k$ form a $\integer[q]$ basis for $QH^*(X)$. A Schubert variety $\Omega_{\lambda}$ determines a class in $QH^*(X)$. It is a result of \cite{B} that this class is still given by $\sigma_{\lambda}$; unlike other flag varieties (see \cite{C-F}, \cite{FGP}), there are no correction terms involving the variable $q$.

It follows that in $QH^*(X)$ there are formulas
\begin{equation}
\sigma_{\lambda}\cdot\sigma_{\mu}=\sum q^mN_{\lambda\mu}^{\nu}(l,k)\sigma_{\nu},
\end{equation}
the sum over $m\geq 0$ and $\nu\subset l\times k$, but with
$ | \nu |  = | \lambda |  + | \mu |  -mn$.
The coefficients $N_{\lambda\mu}^{\nu}(l,k)$ are again nonnegative, this time for 
\lq\lq quantum geometric" reasons: $N_{\lambda\mu}^{\nu}(l,k)$ is the number 
(properly counted) of rational curves of degree $m$ that meet general Schubert varieties
$\Omega_{\lambda}$, $\Omega_{\mu}$, and $\Omega_{\nu^{\vee}}$, where
$\nu^{\vee}=(k-\nu_l,\dots ,k-\nu_1)$ (cf. \cite{B}, \cite{FP}).
The coefficients $N_{\lambda\mu}^{\nu}(l,k)$ are uniquely determined from the presentation (3) of the ring $QH^*(X)$. Similarly, for $0\leq\mu _1,\dots ,\mu_r\leq k$ and $\lambda\subset l\times k$,
\begin{equation}
\sigma_{\mu _1}\cdot\dots\cdot \sigma_{\mu _r}\cdot \sigma_{\lambda}=\sum q^mK_{\lambda\mu}^{\nu}(l,k)\sigma_{\nu},
\end{equation}
the sum over $m\geq 0$, $\nu\subset l\times k$, with
$ | \nu |  = | \lambda |  -mn+\sum_{i=1}^r\mu_i$.
Our goal in this paper is to give explicit formulas for these
{\em quantum Littlewood-Richardson numbers} $N_{\lambda\mu}^{\nu}(l,k)$ and
{\em quantum Kostka numbers} $K_{\lambda\mu}^{\nu}(l,k)$.
Note that they depend on our given numbers $l$ and $k$, as well as on $\lambda$,
$\mu$ and $\nu$. Our result is most satisfying for the quantum Kostka numbers: we show that they are the number of tableaux satisfying a certain condition.
Our algorithm for the quantum Littlewood-Richardson numbers is efficient, but, since it involves signs, it does not show them to be nonnegative.

If $\lambda_1>k$, it follows from the definition that $\sigma_{\lambda}=0$. If
$\lambda_{l+1}>0$, however, $\sigma_{\lambda}$ need not vanish. Our basic algorithm gives a formula for all such polynomials  $\sigma_{\lambda}$. This can be expressed in terms of removing rim $n$-hooks from $\lambda$. There is a {\em rim hook} of ${\lambda}$ corresponding to each of its boxes. It is a {\em rim $n$-hook} if it contains exactly $n$ boxes (which is the same number of boxes as the hook of the corresponding box).

The figure below shows a rim hook with $n=10$.

\vskip 20pt
\epsfxsize 35mm
\centerline{\epsfbox{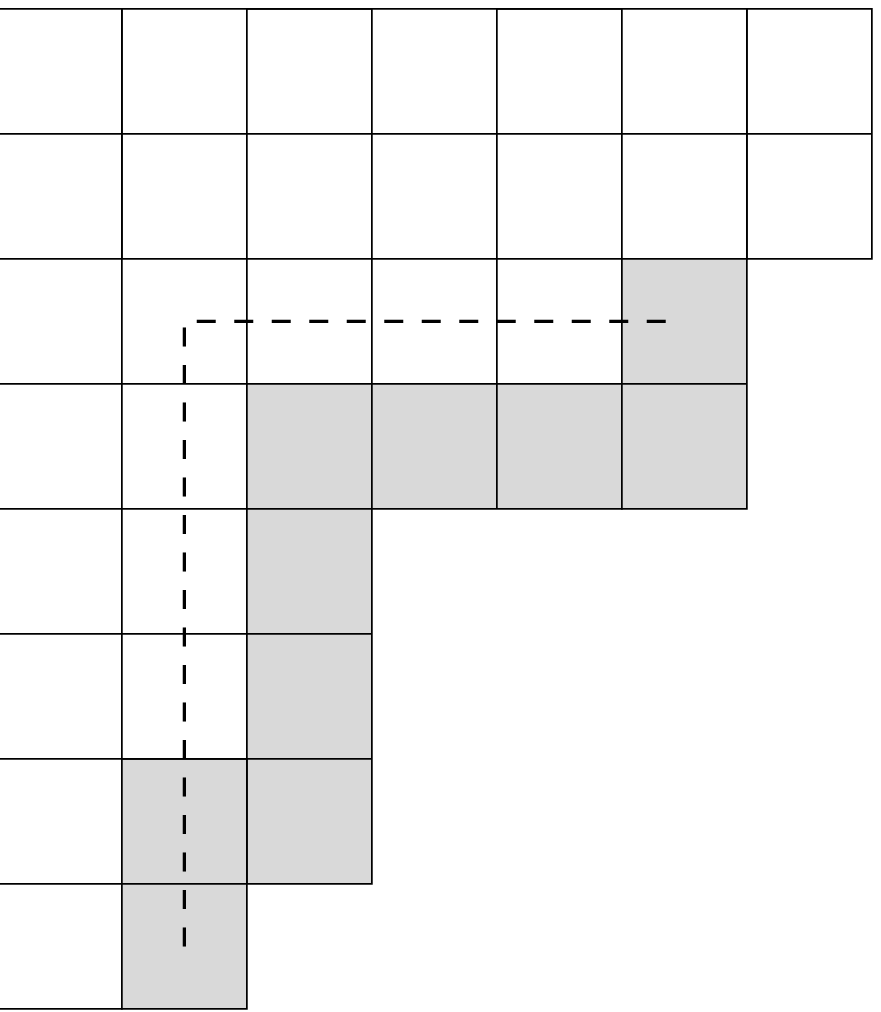}}
\vskip 20pt

Although quantum cohomology is not functorial, isomorphic varieties do have isomorphic quantum cohomology rings. The natural isomorphism of $QH^*(Gr(l,n))$ with
$QH^*(Gr(k,n))$ leads to interesting dual versions of our results, which are discussed in Section 4.

The content of this paper is algebraic and combinatorial; quantum cohomology is
used only for motivation, and for the nonnegativity of the quantum Littlewood-Richardson numbers.

\vskip 10pt

This work was carried out in the stimulating atmosphere of the Mittag-Leffler Institute. We thank F. Sottile and C. Greene for help with the literature on rim hooks. The first author was partially supported by the National Science Foundation. The second author was supported by a Mittag-Leffler Institute postdoctoral fellowship. The third author was supported by a Tage Erlander Guest Professorship, and the National Science Foundation.

\section{The rim hook algorithm}

For positive integers $k,l$, set $n=k+l$, and set
$$\Lambda(l,k)=\integer [q,\sigma_1,\ldots,\sigma_k]/(Y_{l+1}(\sigma),\ldots,Y_{n-1}(\sigma),Y_n(\sigma)+(-1)^kq)\;\; ,$$
where $Y_p(\sigma)={\rm det}(\sigma_{1+j-i})_{1\leq i,j\leq p}$. For any partition $\lambda=(\lambda_1,\ldots,\lambda_r)$, let $\sigma_{\lambda}$ be the image of ${\rm det}(\sigma_{\lambda_i+j-i})_{1\leq i,j\leq r}$ in $\Lambda(l,k)$.

Let $\lambda$ be a partition with $\lambda_1\leq k$, and let $\widetilde{\lambda}_i$ be the number of boxes in the $i^{\rm th}$ column of $\lambda$; i.e., the conjugate partition $\widetilde{\lambda}$ is $(\widetilde{\lambda}_1,\dots ,\widetilde{\lambda}_k)$. One can start at the bottom of any column and move upward and to the right until one has counted off $n$ boxes along the rim. This will be a rim $n$-hook, unless this process ends in a box directly to the left of the last box in some column. We call this an {\em illegal $n$-rim}.

\vskip 20pt
\epsfxsize 25mm
\centerline{\epsfbox{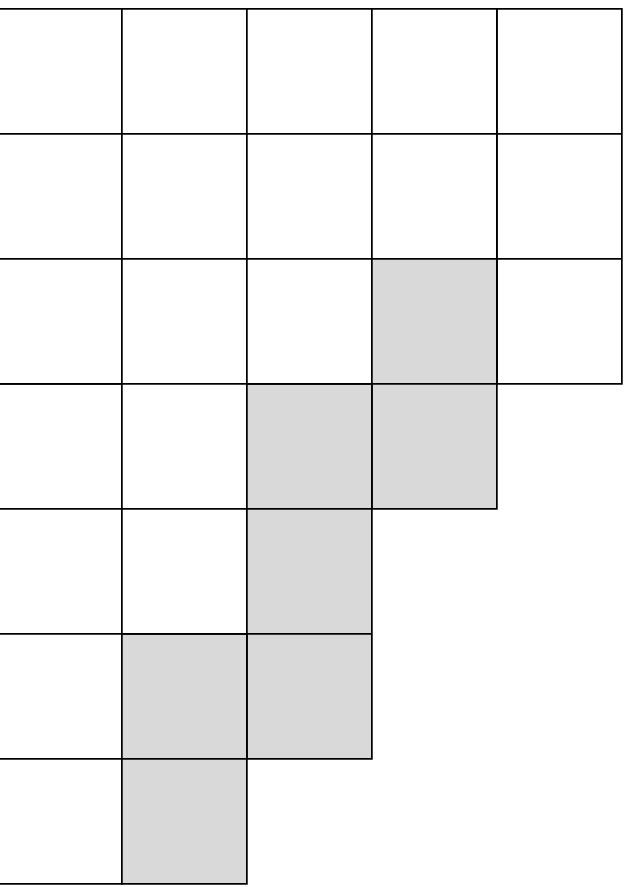}}
\vskip 20pt

If this process starts in column $r$ and ends in column $s$, the column lengths of the shape that remains are
$$(\widetilde{\lambda}_1,\dots ,\widetilde{\lambda}_{r-1},\widetilde{\lambda}_{r+1}-1,
\widetilde{\lambda}_{r+2}-1,\dots ,\widetilde{\lambda}_{s}-1,\widetilde{\lambda}_{r}-r+s-n,
\widetilde{\lambda}_{s+1},\dots ,\widetilde{\lambda}_k).$$

The $n$-rim is illegal exactly when $\widetilde{\lambda}_{r}-r+s-n=\widetilde{\lambda}_{s+1}-1$, or
\begin{equation}
\widetilde{\lambda}_{r}-r-n=\widetilde{\lambda}_{s+1}-(s+1).
\end{equation}
Otherwise $s$ is the unique integer $\geq r$ with
\begin{equation}
\widetilde{\lambda}_{s}-s>\widetilde{\lambda}_{r}-r-n>\widetilde{\lambda}_{s+1}-(s+1).
\end{equation}
The partition $\mu$ obtained after removing such an $n$-rim hook from $\lambda$ is characterized by
\begin{equation}
\widetilde{\mu}=(\widetilde{\lambda}_1,\dots ,\widetilde{\lambda}_{r-1},\widetilde{\lambda}_{r+1}-1,
\dots ,\widetilde{\lambda}_{s}-1,\widetilde{\lambda}_{r}-r+s-n,
\widetilde{\lambda}_{s+1},\dots ,\widetilde{\lambda}_k).
\end{equation}

The {\em width} of the rim $n$-hook is the number $s-r+1$ of columns it occupies. It is also possible that there is no $n$-rim starting in column $r$, which happens if $\widetilde{\lambda}_r+\lambda_1-r<n$; so if
$\widetilde{\lambda}_1+\lambda_1\leq n$, then $\lambda$ contains no $n$-rim.

\begin{mlemma}{\nonumber} $({\rm A})$ If $\lambda$ contains an illegal $n$-rim, or if $\lambda_{l+1}>0$ and $\lambda$ contains no $n$-rim, then $\sigma_{\lambda}=0$.

$({\rm B})$ If $\mu$ is the result of removing a rim $n$-hook from $\lambda$, then
$$\sigma_{\lambda}=(-1)^{k-w}q\sigma_{\mu},$$
where $w$ is the width of the rim $n$-hook removed.
\end{mlemma}

\prf Let $y_r=\sigma_{(1^r)}$ for all integers $r\geq 1$, with $y_0=1$ and $y_i=0$ for $i<0$. We know by (1) and (3) that $y_i=0$ for $l<i<n$, and $y_n=(-1)^{k-1}q$. By expanding the determinant expression for $y_r$ along the top row, we have
$$y_r=\sigma_1y_{r-1}-\sigma_2y_{r-2}+\dots +(-1)^{k-1}\sigma_ky_{r-k}.$$
From this it follows by induction that
\begin{equation}
y_{mn+j}=(-1)^{m(k-1)}q^my_j \qquad {\rm for} \quad 0\leq j\leq n-1,\ m\geq 0.
\end{equation}
For any sequence $m=(m_1,\dots ,m_k)$ of integers, set
\begin{equation}
\tau_m={\rm det}(y_{m_i+j-i})_{1\leq i,j\leq k}.
\end{equation}
If $\lambda$ is a partition with $\lambda_1\leq k$, it is a basic identity for
symmetric polynomials (cf. \cite{M}, $(2.9')$) that
\begin{equation}
\sigma_{\lambda}=\tau_{\widetilde{\lambda}}
\end{equation}
If $1\leq r<s\leq k$ and the $r^{\rm th}$ row of the matrix $(y_{m_i+j-i})_{1\leq i,j\leq k}$ is successively interchanged with the $(r+1)^{\rm st}$ row, then the $(r+2)^{\rm nd}$ row, and so on to the $s^{\rm th}$ row, one finds
\begin{equation}
\tau_m=(-1)^{s-r}\tau_{m'}\ ,
\end{equation}
where $m'=(m_1,\dots ,m_{r-1},m_{r+1}-1,\dots ,m_s-1,m_r-r+s,m_{s+1},\dots ,m_k)$.

Suppose the $n$-rim starting in column $r$ of $\lambda$ ends in column $s$. Set
$m=\widetilde{\lambda}$. If this rim is illegal, it follows from (11) that
$(y_{m_i'+j-i})$ has two identical rows, so $\sigma_{\lambda}=\tau_m=0$. This proves the first part of $({\rm A})$. If $\widetilde{\lambda}_r\geq l$ for some $r$ between $1$ and $k$, applying (14) to the
entries of row $s$ of $(y_{m_i'+j-i})$ yields
\begin{equation}\tau_{m'}=(-1)^{k-1}q\tau_{m''}\ ,
\end{equation}
where $m''=(m'_1,\dots ,m'_{s-1},m'_s-n,m'_{s+1},\dots ,m'_k)$. If the $n$-rim starting in column $r$ of $\lambda$ is a rim hook, and it ends in column $s$,
then by (17) and (18)
$$\tau_m=(-1)^{k-1+s-r}q\tau_{m'''}\ ,$$
where $m'''=(m_1,\dots ,m_{r-1},m_{r+1}-1,\dots ,m_s-1,m_r-r+s-n,m_{s+1},\dots ,m_k)$. Part $({\rm B})$ of the lemma follows from this, using (13) and (16).

On the other hand, if $\lambda_{l+1}>0$ and $\lambda$ contains no $n$-rim, we have $m_1>l$, $m_1+\lambda_1\leq n$, with $m=\widetilde{\lambda}$. Then $\sigma_{\lambda}={\rm det}(y_{m_i+j-i})_{1\leq i,j\leq \lambda_1}$, which vanishes since its top row is zero. This completes the proof of $({\rm A})$.\qed

\vskip 12pt

Given any partition $\lambda$ with $\lambda_1\leq k$, one can succesively apply the lemma, removing rim $n$-hooks until one arrives at a partition $\mu$ which has no rim $n$-hooks. If $\mu\subset l\times k$, and $R_1,\dots ,R_m$ are the rim hooks removed, then
\begin{equation}
\sigma_{\lambda}=\epsilon(\lambda/\mu)q^m\sigma_{\mu},
\end{equation}
where $\epsilon(\lambda/\mu)=(-1)^{\sum(k-{\rm width}(R_i))}$. If $\mu$ is not contained in the $l\times k$ rectangle, then
\begin{equation}
\sigma_{\lambda}=0.
\end{equation}

One sees from (19) that if $\mu\subset l\times k$, then $\mu$ is independent of choice of the $m$ rim hooks removed, as is $\sum{\rm width}(R_i)$ mod($2$). This is a general fact: $\mu$ is the {\em $n$-core} of $\lambda$, which is uniquely determined by $\lambda$ and $n$ (\cite{JK}, \S 2.7, cf. \cite{M}, \S 1).

Consider the specialization of (6) by $s_{\lambda}\longmapsto \sigma_{\lambda}$,
using (19) and (20). This gives
\begin{equation}
\sigma_{\lambda}\cdot\sigma_{\mu}=\sum_{m,\nu}q^mN_{\lambda\mu}^{\nu}(l,k)\sigma_{\nu},\;\;\; 
N_{\lambda\mu}^{\nu}(l,k)=\sum_{\varrho}\epsilon(\varrho/\nu)N_{\lambda\mu}^{\varrho}.
\end{equation} 
Here $\lambda,\mu\subset l\times k$, the first sum is over $m\geq 0$ and $\nu\subset l\times k$ with $$|\lambda|+|\mu|=|\nu|+mn,$$ 
the second sum is over $\varrho$ obtained from $\nu$ by adding $m$ rim $n$-hooks,  
and $N_{\lambda\mu}^{\varrho}$ are the {\em classical} Littlewood-Richardson coefficients. This gives a formula for the quantum Littlewood-Richardson numbers in terms of the classical numbers. From the identification of $\Lambda(l,k)$ with $QH^*(X)$ we deduce:

\begin{corollary}
For $\lambda,\mu,\nu\subset l\times k$, with $|\lambda|+|\mu|=|\nu|+mn$ for some $m\geq 0$, the number of rational curves of degree $m$ that meet general translates of the Schubert varieties
$\Omega_{\lambda}$, $\Omega_{\mu}$, and $\Omega_{\nu^{\vee}}$ is equal to
$$\sum\epsilon(\varrho/\nu)N_{\lambda\mu}^{\varrho},$$
where the sum is over all $\varrho$ with $\varrho_1\leq k$ that can be obtained from $\nu$ by adding $m$ rim $n$-hooks. In particular, $\sum\epsilon(\varrho/\nu)N_{\lambda\mu}^{\varrho}\geq 0$.\qed
\end{corollary}

\begin{numexample}
{\rm $k=5$, $l=5$, $\lambda=(5,4,4,2,2)$, $\mu=(3,2,1)$, $\nu=(2,1)$:}

\vskip 20pt
\epsfxsize 13cm
\centerline{\epsfbox{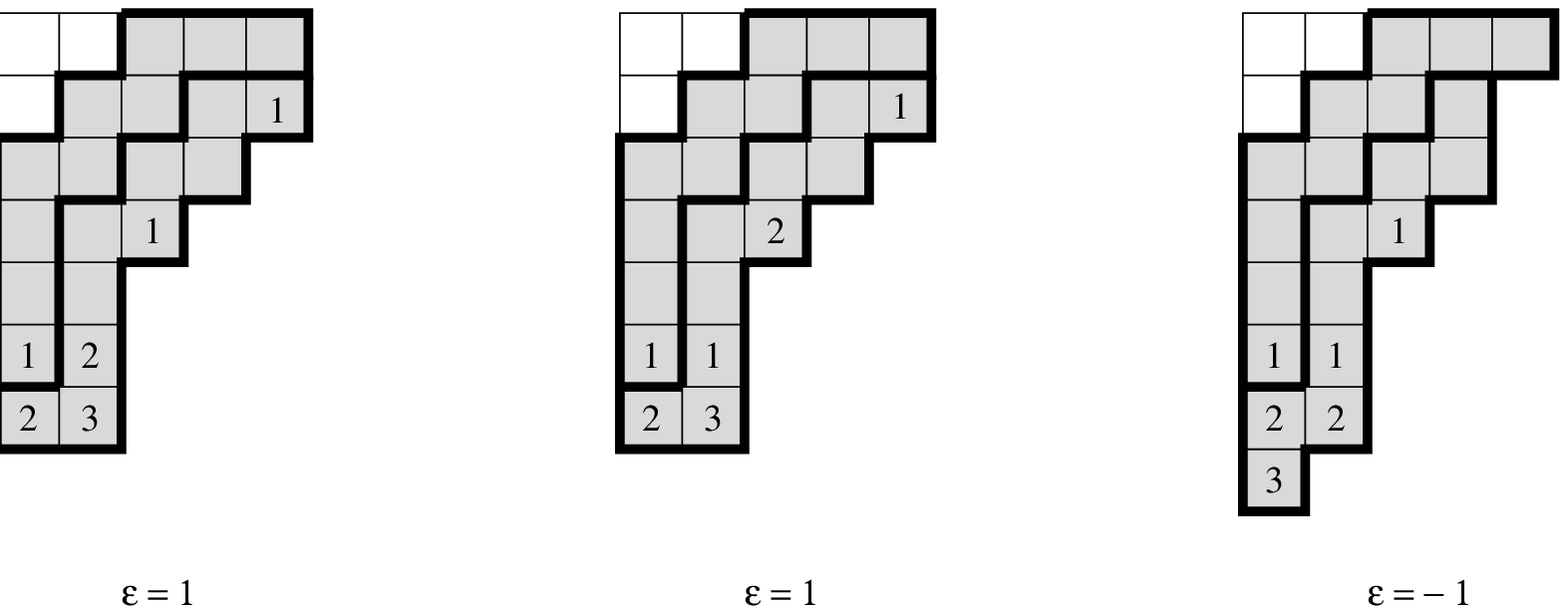}}
\vskip 20pt

\centerline{$N_{\lambda\mu}^{\nu}(5,5)=2-1=1$.}
\end{numexample}

\begin{numexample}
{\rm This example exhibits the dependence of $N_{\lambda\mu}^{\nu}(l,k)$ on both $l$ and $k$, as well as an instance where a quantum Littlewood-Richardson number
(which is a nontrivial sum of signed classical ones) vanishes.
 Let $\lambda=(3,3,2,1)$, $\mu=(4,3,2,1)$, $\nu=(4,2,2,1)$. By removing a rim $10$-hook from each of the partitions $\varrho_1=(6,5,3,3,2)$, $\varrho_2=(5,5,3,3,2,1)$, and $\varrho_3=(4,4,3,3,2,1,1,1)$, one obtains the shape $\nu$, and these are the only partitions with this property that appear in the product $s_{\lambda}\cdot s_{\mu}$ . The classical 
Littlewood-Richardson coefficients are
$$N_{\lambda\mu}^{\varrho_1}=6,\qquad N_{\lambda\mu}^{\varrho_2}=8,\qquad N_{\lambda\mu}^{\varrho_3}=2.$$
Our algorithm gives
$$
N_{\lambda\mu}^{\nu}(4,6)=N_{\lambda\mu}^{\varrho_1}-N_{\lambda\mu}^{\varrho_2}+N_{\lambda\mu}^{\varrho_3}=6-8+2=0,
$$
$$
N_{\lambda\mu}^{\nu}(5,5)=N_{\lambda\mu}^{\varrho_2}-N_{\lambda\mu}^{\varrho_3}=
8-2=6,
$$
and
$$
N_{\lambda\mu}^{\nu}(6,4)=N_{\lambda\mu}^{\varrho_3}=2.
$$
}
\end{numexample}

\vskip 10pt

The quantum Pieri formula, for multiplying any $\sigma_{\lambda}$ by a class
$\sigma_p$, was given in \cite{B}: for $\lambda\subset l\times k$, $1\leq p\leq k$,
\begin{equation}
\sigma_p\cdot\sigma_{\lambda}=\sum\sigma_{\mu}+\sum q\sigma_{\nu}.
\end{equation}
Here the first (classical) sum is over $\mu$ with
$$k\geq \mu_1\geq\lambda_1\geq\mu_2\geq\lambda_2\geq\dots\geq\mu_l\geq\lambda_l\geq 0,$$
and $|\mu|=|\lambda|+p $; the second sum is over $\nu$ with
$$\lambda_1-1\geq\nu_1\geq\lambda_2-1\geq\nu_2\geq\dots\geq\lambda_l-1\geq\nu_l\geq 0,$$ and $|\nu|=|\lambda|+p-n$. These $\nu$ are the partitions obtained from $\lambda$ by removing $n-p$ boxes from its border rim, taking at least one box from each of the $l$ rows.

The proof in \cite{B} combined algebra with some quantum geometry. This Pieri formula follows easily from the Main Lemma. Indeed, one has the formula for multiplying Schur functions
$$s_p\cdot s_{\lambda}=\sum s_{\mu},$$
the sum over $\mu$ obtained from $\lambda$ by adding $p$ boxes, no two in the same column. Terms with $\mu_1>k$, or with $\mu_1<k$ and $\mu_{l+1}>0$ vanish by $({\rm A})$ of the Main Lemma. Those with $\mu_{l+1}=0$ give the first sum in (22). Those with $\mu_{l+1}>0$ and $\mu_1=k$ have a rim $n$-hook that occupies all $k$ columns; by $({\rm B})$ of the Main Lemma, these give the second sum in (22).
\begin{numexample}
{\rm $k=4$, $l=4$, $\lambda=(2,2,1,1)$, $p=3$:}

\vskip 20pt
\epsfxsize 40mm
\centerline{\epsfbox{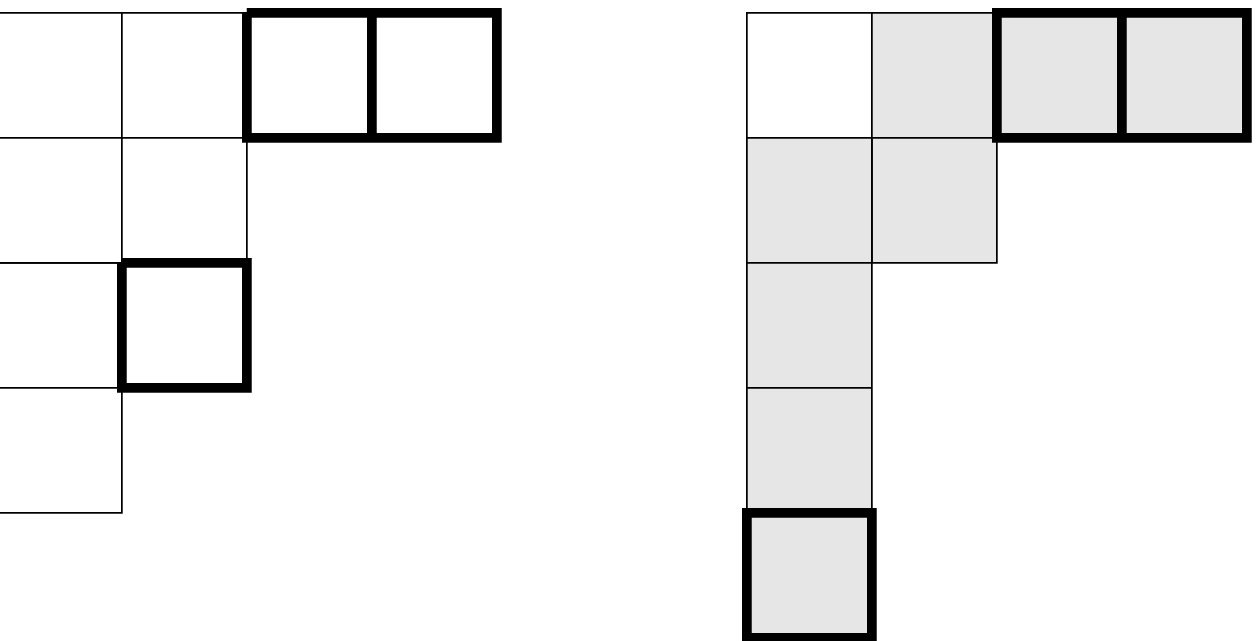}}
\vskip 20pt

\centerline{$\sigma_3\cdot\sigma_{(2,2,1,1)}=\sigma_{(4,2,2,1)}+q\sigma_1$.}
\end{numexample}

\bigskip

\section{Quantum Kostka numbers}

\bigskip

For a partition $\nu\subset l\times k$, and a nonnegative integer $m$, we denote by $\nu[m]$ the partition obtained from $\nu$ by adding $m$ rim $n$-hooks to $\nu$, each starting in the first column and ending in the $k^{\rm th}$ column. Note that the height of each added rim is $l+1$, i.e., its lowest box in the first column is $l$ rows below its highest box in the $k^{\rm th}$ column.

Let $\lambda\subset\varrho$ be partitions, with $\varrho_1\leq k$. We call a (semistandard) tableau $T$ on $\varrho/\lambda$ {\em proper} if, for each entry $i$ of $T$ in the first column and row $l+p$, $p\geq 1$, the box in the 
$k^{\rm th}$ column and row $p$ is either in $\lambda$, or it is in $\varrho/\lambda$ and it has an entry of $T$, at most equal to $i$.

\vskip 20pt
\epsfxsize 9cm
\centerline{\epsfbox{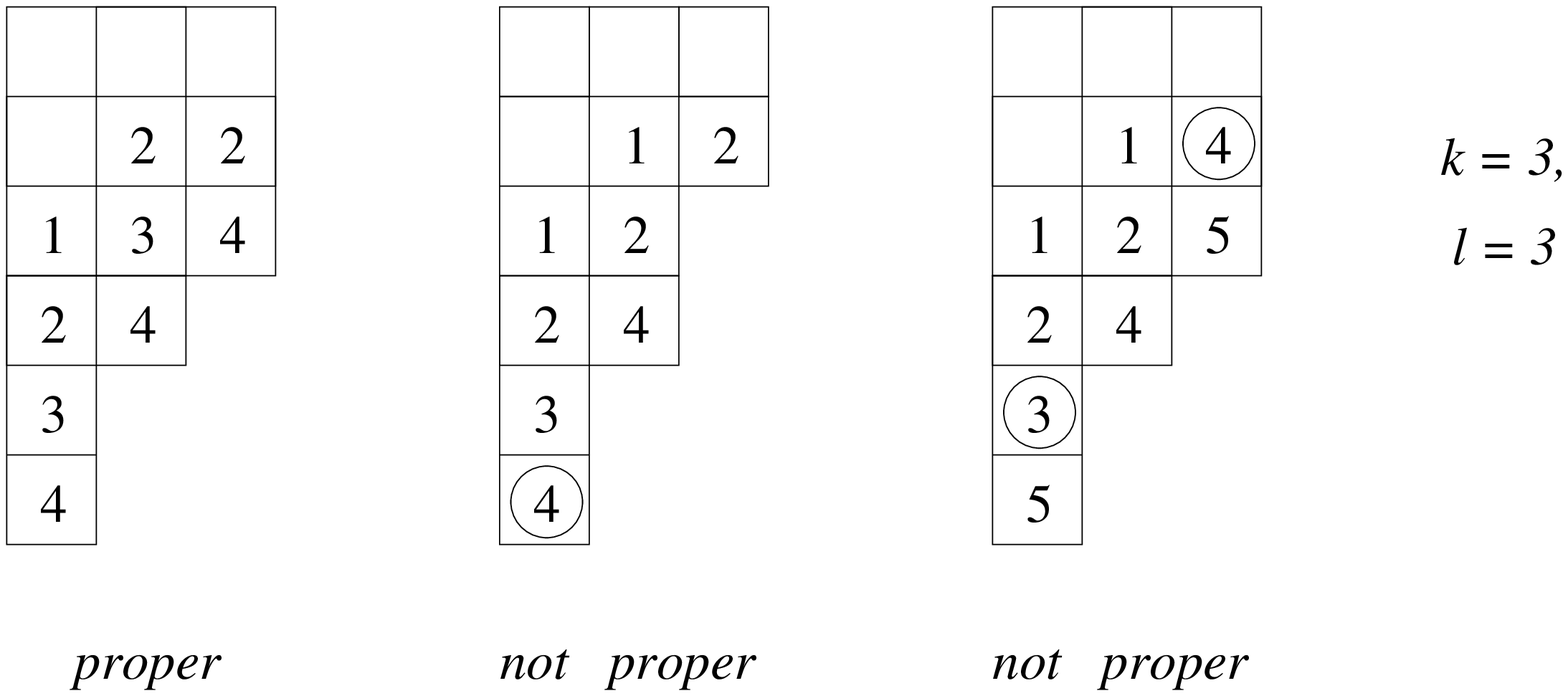}}
\vskip 20pt

For any sequence $\mu=(\mu_1,\dots ,\mu_r)$ of nonnegative integers, each at most $k$, and partitions $\lambda,\nu\subset l\times k$, with $|\lambda|+\sum\mu_i=|\nu|+mn$ for some $m\geq 0$, define the {\em quantum Kostka number} $K_{\lambda\mu}^{\nu}(l,k)$ to be the number of proper tableaux on $\nu[m]/\lambda$ with content $\mu$. When $m=0$, this is the 
classical Kostka number $K_{\lambda\mu}^{\nu}$ (the coefficient of
$s_{\nu}$ in $s_{\mu_r}\cdot\ldots\cdot s_{\mu_1}\cdot s_{\lambda}$).

\begin{proposition} For $\lambda\subset l\times k$ and $0\leq \mu_1,\dots ,\mu_r\leq k$,
$$\sigma_{\mu_r}\cdot\dots\cdot\sigma_{\mu_1}\cdot\sigma_{\lambda}=\sum
K_{\lambda\mu}^{\nu}(l,k)q^m\sigma_{\nu},$$
the sum over $\nu\subset l\times k$ and $m\geq 0$, with
$|\nu|+mn=|\lambda|+\sum\mu_i$.
\end{proposition}

Note in particular that $K_{\lambda\mu}^{\nu}(l,k)$ is independent of the ordering of the terms in $\mu=(\mu_1,\dots ,\mu_r)$.

The proof will follow from three lemmas. Fix $\lambda$ and $\mu=(\mu_1,\dots ,\mu_r)$ as in the proposition. For any tableau $T$ of content $\mu$ on a skew shape $\varrho/\lambda $, and $1\leq i\leq r$, let $T_i$ be the subtableau of $T$ obtained by discarding all entries strictly larger than $i$; $T_i$ is a tableau of content $(\mu_1,\dots ,\mu_i)$ on some shape $\varrho(i)/\lambda$. 
Set $\varrho(0)=\lambda$.
\begin{lemma}
A tableau $T$ of content $\mu$ on $\varrho/\lambda$ is proper if and only if $$\widetilde{\varrho(i)}_1-\widetilde{\varrho(i)}_k\leq l,\ \ \ {\it for}\ \ 1\leq i\leq r.$$
\end{lemma}

\prf The condition of the definition is satisfied for an entry $i$ in the first column of $T$ if and only if the $k^{\rm th}$ column of $\varrho(i)$ is no more than $l$ boxes shorter than its first column.\qed

\vskip 12pt

For any $T$ that is not proper, let $U(T)=T_i$, where $i$ is minimal not satisfying the condition of the lemma. We call $U$ the {\em improper kernel} of $T$, and $T$ a {\em descendent} of $U$ if $U(T)=U$. Note that each $\varrho(i)$ has at most one more box in any column than $\varrho(i-1)$. It follows that if $U=U(T)$ has shape $\alpha/\lambda$, then (since $\alpha=\varrho(i)$ and $i$ is minimal)
\begin{equation}
\widetilde{\alpha}_1-\widetilde{\alpha}_k=l+1.
\end{equation}

For any tableau $T$ on a shape $\varrho/\lambda$, write $\sigma_T$ for the element
$\sigma_{\varrho}$ in $\Lambda(l,k)$.

\begin{lemma}
Let $U$ be a tableau on $\alpha/\lambda$ with content $(\mu_1\dots ,\mu_i)$ for some $1\leq i\leq r$. Suppose $\widetilde{\alpha}_1-\widetilde{\alpha}_k=l+1$.
Then $$\sum\sigma_T=0,$$
the sum over all $T$ of content $\mu$ with $T_i=U$.
\end{lemma}

\prf The sum $\sum\sigma_T$ in the statement is equal to
$$\sigma_{\mu_r}\cdot\sigma_{\mu_{r-1}}\cdot\dots\cdot\sigma_{\mu_{i+1}}\cdot\sigma_{\alpha}$$ in $\Lambda(l,k)$. It therefore suffices to show that $\sigma_{\alpha}=0$.
From the assumption $\widetilde{\alpha}_1-\widetilde{\alpha}_k=l+1$ it follows that the $n$-rim starting at the bottom of the first column will end in the $(k-1)^{\rm st}$ column, just to the left of the last box in the $k^{\rm th}$ column.

\vskip 20pt
\epsfxsize 8cm
\centerline{\epsfbox{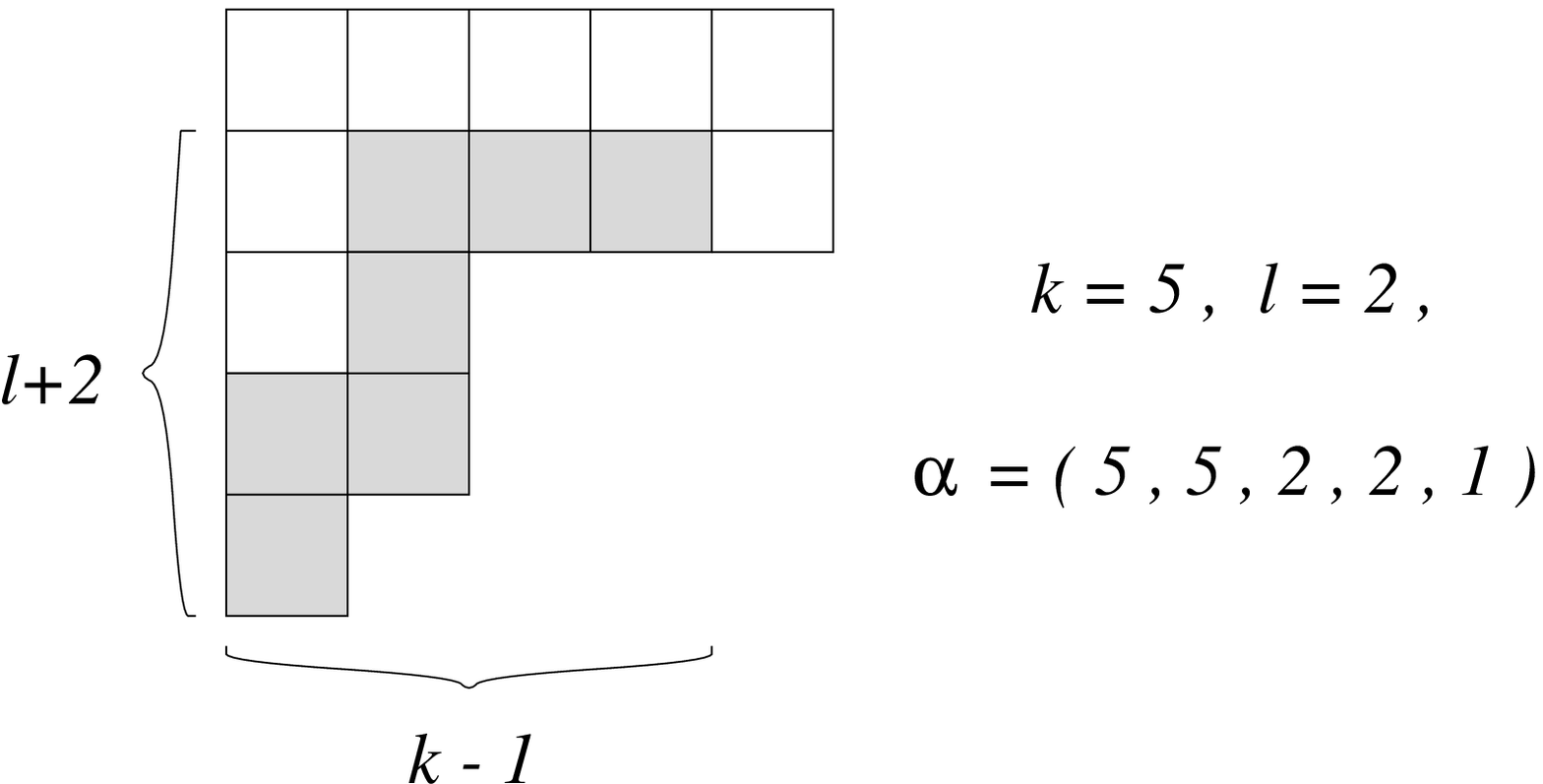}}
\vskip 20pt

\noindent This is an illegal $n$-rim, so $\sigma_{\alpha}=0$ by $({\rm A})$ of the Main Lemma.\qed

\begin{lemma}
Let $\varrho$ be a partition with $\varrho_1\leq k$ and $\widetilde{\varrho}_1-\widetilde{\varrho}_k\leq l$. Then there is a unique partition $\nu\subset l\times k$ and $m\geq 0$ with $\varrho=\nu[m]$. Moreover, $\sigma_{\varrho}=q^m\sigma_{\nu}$ in $\Lambda(l,k)$.	
\end{lemma}

\prf If $\widetilde{\varrho}_1\leq l$, then $\nu=\varrho$ and $m=0$. Otherwise, the condition $\widetilde{\varrho}_1-\widetilde{\varrho}_k\leq l$ implies that $\varrho$ contains a rim $n$-hook of width $k$ and height $l+1$. Removing this rim hook leaves a shape $\beta$ with $\widetilde{\beta}_1-\widetilde{\beta}_k\leq l$, and
$\sigma_{\varrho}=q\sigma_{\beta}$ by $({\rm B})$ of the Main Lemma. The proof concludes by induction on the number of boxes.\qed

\vskip 12pt

To prove the proposition, from the classical expansion of Schur functions (cf. \cite{F2}, \S 2, or \cite{M}, \S I.5), we have
$$\sigma_{\mu_r}\cdot\dots\cdot\sigma_{\mu_1}\cdot\sigma_{\lambda}=\sum
\sigma_{T},$$
the sum over tableaux $T$ of content $\mu$ on some shape $\varrho(T)/\lambda$. The shape $\varrho=\varrho(T)$ for any proper $T$ satisfies the condition of Lemma 3 (by Lemma 1), so the sum of the $\sigma_T$ for $T$ proper gives the right side of the equation in the proposition. The $T$ that are not proper are divided into equivalence classes indexed by their improper kernels $U=U(T)$. Each sum $\sum\sigma_T$ over $T$ descended from $U$ vanishes by Lemma 2, which concludes the proof of the proposition.\qed

\vskip 12pt

This proposition gives another algorithm for computing the quantum product
$\sigma_{\lambda}\cdot\sigma_{\mu}$, for $\lambda$ and $\mu$ partitions in $l\times k$:
Expand
$$\sigma_{\mu}={\rm det}(\sigma_{\mu_i+j-i})=\sum_{\tau\in S_l}{\rm sign}(\tau)\prod_{i=1}^l\sigma_{\mu_i+\tau(i)-i}$$ and use the proposition to expand each $\prod\sigma_{\mu_i+\tau(i)-i}\cdot\sigma_{\lambda}$. This writes $\sigma_{\lambda}\cdot\sigma_{\mu}$ as an alternating sum of sums of $\sigma_T$'s as $T$ runs over proper tableaux of content $(\mu_1+\tau(1)-1,\dots ,\mu_l+\tau(l)-l)$ on shapes $\nu[m]/\lambda$. Note that among the sum with $\tau=id$ are all the tableaux whose word is a reverse lattice word. It is known from the classical Littlewood-Richardson rule that the other tableaux cancel in such a sum, but it is unclear how proper and nonproper tableaux behave under
such a cancellation.

\bigskip

\section{Duality}

\bigskip

Although the functoriality of ordinary cohomology does not extend to quantum
cohomology, isomorphic varieties do have isomorphic quantum cohomology rings. The canonical isomorphism of $Gr(l,n)$ with $Gr(k,n)$ gives an isomorphism of $QH^*(Gr(l,n))$ with $QH^*(Gr(k,n))$, which takes the class of a Schubert variety $\Omega_{\lambda}$ to the class of a Schubert variety $\Omega_{\widetilde{\lambda}}$, for $\lambda\subset l\times k$. Our basic algorithm from Section 2 is not at all invariant under the involution $\lambda\mapsto\widetilde{\lambda}$, however, so this isomorphism leads to some interesting combinatorial identities.

We begin with an algebraic proof of this isomorphism. To avoid cofusion, set
$$\Lambda(k,l)=\integer[q,\tau_1,\ldots,\tau_l ]/(Y_{k+1}(\tau),\ldots Y_{n-1}(\tau), Y_n(\tau)+(-1)^lq).$$ 

\begin{numproposition}
The map $\sigma_i\mapsto Y_i(\tau)$, $1\leq i\leq k$, determines an isomorphism $\Lambda(l,k)\stackrel{\sim}{\rightarrow}\Lambda(k,l)$, with inverse taking $\tau_i$ to $Y_i(\sigma)$, $1\leq i\leq l$. For any $\lambda\subset l\times k$,
this isomorphism takes $\sigma_{\lambda}={\rm det}(\sigma_{\lambda_i+j-i})_{1\leq i,j\leq l}$ to $\tau_{\widetilde{\lambda}}={\rm det}(\tau_{\widetilde{\lambda}_i+j-i})_{1\leq i,j\leq k}$.
\end{numproposition}

\prf Consider the ring 
$$\Lambda=\integer [q,\sigma_1,\ldots,\sigma_k,\tau_1,\ldots,\tau_l ]/I\; ,$$
where $I$ is generated by the elements
$$\sigma_1-\tau_1,\;\sigma_2-\sigma_1\tau_1+\tau_2,\;\ldots,\;
\sigma_k\tau_{l-1}-\sigma_{k-1}\tau_l,\;\sigma_k\tau_l-q.$$
These relations identify $\tau_i$ with $Y_i(\sigma)$ for $1\leq i\leq l$, and then prescribe that $Y_i(\sigma)=0$ for $l<i<n$, and that $Y_n(\sigma)=(-1)^{k-1}q$. This identifies $\Lambda$ with $\Lambda(l,k)$. By symmetry, it also identifies $\Lambda$ with $\Lambda(k,l)$, and the composite
$\Lambda(l,k)\cong\Lambda\cong\Lambda(k,l)$ is the isomorphism of the proposition.

The fact that $\sigma_{\lambda}$ and $\tau_{\widetilde{\lambda}}$ correspond under this isomorphism follows from the general identity \cite{M}, $(2.9')$, and the fact that all $\sigma_i$ (resp. $\tau_i$) occuring in the matrix for $\sigma_{\lambda}$ (resp. $\tau_{\widetilde{\lambda}}$) have $i<n$.\qed

\vskip 12pt

The following dual of the quantum Pieri formula follows from this proposition. We include a direct proof, for contrast with the proof of the Pieri formula in Section 2.

\begin{numproposition} Let $p\leq l$ be a positive integer and let $\lambda\subset l\times k$ be a partition.
The following holds in $\Lambda(l,k)$:
\begin{equation}\label{pieri2}
\sigma_{\lambda}\cdot\sigma_{(1^p)}=\sum\sigma_{\mu}+q\sum\sigma_{\nu},
\end{equation}
the first sum over $\mu\subset l\times k$, obtained by adding
$p$ boxes to $\lambda$, with no two in the same row, and the second sum
over $\nu$ obtained from $\lambda$ by removing $n-p$ boxes in its
border rim, with at least one from each of the $k$ columns.
Note that there are no such $\nu$ if $\lambda_1<k$.

Equivalently, the partitions $\nu$ occuring in the second sum can be characterized by $ | \nu |  = | \lambda |  +p-n$ and
$$
\widetilde{\lambda}_1-1\geq\widetilde{\nu}_1\geq\widetilde{\lambda}_2-1\geq\widetilde{\nu}_2
\geq\dots\geq\widetilde{\lambda}_k-1\geq\widetilde{\nu}_k\geq 0.$$ \end{numproposition}

\prf The classical Pieri formula for Schur polynomials states that
\begin{equation}\label{star}
s_{\lambda}\cdot s_{(1^p)}=\sum s_{\pi},
\end{equation}
sum over all $\pi$ obtained by adding $p$ boxes to $\lambda$, with no two in the same row. We need to show that the terms with $\pi_{l+1}>0$ give the second sum in (\ref{pieri2}).

Let $\cP$ denote the set of partitions $\pi$ that occur in (\ref{star}), with $\pi_1\leq k$, $\pi_{l+1}>0$, and the $n$-core of $\pi$
contained in the $l\times k$ rectangle. Only partitions in $\cP$ can contribute to the second sum in (\ref{pieri2}), 
by $({\rm A})$ of the Main Lemma. The formula (\ref{pieri2}) will follow therefore from the following two claims:
\vskip 10pt
\noindent{\bf Claim 1.} If $\nu$ is as in the proposition, then there exists
a unique partition $\pi\in\cP$ such that $\nu$ is obtained from $\pi$ by removing an $n$-rim. Moreover, $\epsilon(\pi/\nu)=1$.

\vskip 12pt

\noindent {\bf Claim 2.} Let $\pi\in\cP$ be any partition not arising in 1. Then there exists exactly one other partition $\pi '\in\cP$
with the same $n$-core $\varrho$ as $\pi$. Moreover, $\epsilon(\pi/\varrho)=
-\epsilon(\pi '/\varrho)$.

\vskip 10pt

{\it Proof of Claim $1$:} The partition $\pi$ is obtained by adding $p$ boxes
to $\lambda$ as follows:

\noindent For each $i\in\{ 2,3,\dots ,k\}$, add $\widetilde{\nu}_{i-1}-\widetilde{\lambda}_i+1$ boxes to the $i^{\rm th}$ column of $\lambda$. The condition $\widetilde{\lambda}_1-1\geq\widetilde{\nu}_1\geq\widetilde{\lambda}_2-1\geq\widetilde{\nu}_2
\geq\dots\geq\widetilde{\lambda}_k-1\geq\widetilde{\nu}_k\geq 0$ ensures that no two of these new boxes are in the same row; the total number of boxes added according to this recipe is 
$$\sum_{i=2}^k\widetilde{\nu}_{i-1}-\widetilde{\lambda}_i+1= | \nu |  -\widetilde{\nu}_k-
 | \lambda |  +\widetilde{\lambda}_1+k-1=p-(l+1-\widetilde{\lambda}_1+\widetilde{\nu}_k).$$ Now add $l+1-\widetilde{\lambda}_1+\widetilde{\nu}_k$ boxes to the first column
of $\lambda$ to get a partition $\pi$.

Then $\pi$ contains the rim $n$-hook consisting of these $p$ added boxes together with the $n-p$ boxes in $\lambda/\nu$; its width is $k$, so $\epsilon(\pi/\nu)=1$. This is the only way a $\pi$ in $\cP$ can arise from $\nu$ by adding a rim $n$-hook.

\vskip 12pt

{\it Proof of Claim $2$:} Let $\pi$ be a partition in $\cP$. The $n$-rim of $\pi$ starting at the bottom of the first column is a rim $n$-hook whose removal gives the $n$-core $\varrho$ of $\pi$. This rim hook consists of some boxes in the rim of $\lambda$
and some of the $p$ added boxes.
Let $t$ be the maximum index such that this rim hook contains at least one box from the $t^{\rm th}$ column of $\lambda$. 

If $t=k$, then $\pi$ is one of the partitions considered in 1. 
We therefore assume that $t\neq k$. The rim hook must either $(a)$ end with a box from $\lambda$ in the 
$t^{\rm th}$ column of $\pi$, or $(b)$ end with an added box in the 
$(t+1)^{\rm st}$ column. 

In case $(a)$, $\pi '$ is obtained from $\pi$ by removing
$\widetilde{\varrho}_t-\widetilde{\pi}_{t+1}+1$ boxes from the first column, and adding them to the $(t+1)^{\rm st}$ column of $\pi$; note that
$\pi '$ satisfies $(b)$, has the same $n$-core $\varrho$, and $\epsilon(\pi/\varrho)=
-\epsilon(\pi '/\varrho)$.

Conversely, if $\pi$ satisfies $(b)$, the only $\pi '$ with the same $n$-core
as $\pi$ is obtained from $\pi$ by removing $\widetilde{\pi}_{t+1}-\widetilde{\varrho}_{t+1}$ boxes from the 
$(t+1)^{\rm st}$ column and adding them to the first column.  \qed

\vskip 12pt

We turn next to a dual form of the proposition in Section 3; this is {\em not} obtained by applying the duality isomorphism of Proposition 4.1.
Call a filling $T$ of the boxes of a skew shape $\varrho / \lambda$
a {\em conjugate tableau} if its entries are strictly increasing 
across rows and weakly increasing down columns. Call a conjugate
tableau {\em proper} if for each entry $i$ of $T$ in the first column and row $l+p$, $p \geq 1$, the box in the $k^{\rm th}$ column and row $p$ is either in $\lambda$, or it has an entry of $T$ {\em strictly}
smaller than $i$. 

For a sequence $\mu =(\mu_1,\ldots ,\mu_r)$ of nonnegative integers,
and partitions $\lambda ,\nu \subset l \times k$, with
$|\lambda|+\sum\mu_i=|\nu|+mn$, $m\geq 0$, define the {\em conjugate 
quantum Kostka number} $\widetilde{K}_{\lambda\mu}^{\nu}(l,k)$ to be the
number of proper conjugate tableaux on $\nu[m] / \lambda$ with
content $\mu$.

\begin{numproposition}
For $\lambda \subset l \times k$ and $0\leq \mu_1\leq \ldots \leq \mu_r\leq l$,
$$
\sigma_{(1^{\mu_r})}\cdot\ldots\cdot\sigma_{(1^{\mu_1})}\cdot\sigma_{\lambda}=
\sum\widetilde{K}_{\lambda\mu}^{\nu}(l,k) q^m \sigma_{\nu},
$$
the sum over $\nu\subset l \times k$ and $m\geq 0$ with
$|\nu|+mn=|\lambda|+\sum\mu_i$.
\end{numproposition}

\prf As in Section 3, for any conjugate tableau $T$ on shape $\varrho/\lambda$ with content $\mu$, let $T_i$ be the part of $T$ of content $(\mu_1,\ldots,\mu_i)$, of shape $\varrho(i)/\lambda$.
This time $T$ is proper if and only if, for $1\leq i\leq r$, $\widetilde{\varrho(i)}_1-\widetilde{\varrho(i)}_k\leq l$, with the inequality strict if $T$ has an entry $i$ in the $k^{\rm th}$ column.

If $T$ is not proper, let $U(T)=T_i$ for $i$ minimal not satisfying the condition of the lemma; call $T$ a descendent of $U(T)$. We consider conjugate tableaux $U$ of content $(\mu_1,\ldots,\mu_i)$ which are not proper, but such that $U_{i-1}$ (of content $(\mu_1,\ldots,\mu_{i-1})$) is proper and fixed.
Consider the $n$-rim of the shape $\varrho(U)$ that starts in the first column. If this rim is illegal, then $\sigma_{\varrho(U)}=0$. It follows as in Lemma 2 of \S 3 that $\sum\sigma_T=0$, the sum over all $T$ descended from $U$.

For each $U$ whose $n$-rim is a rim hook we will construct another conjugate tableau $U'$ of content $(\mu_1,\ldots,\mu_i)$, with $U'_{i-1}=U_{i-1}$, whose $n$-rim is also a rim hook. The results of removing the rim $n$-hooks from $U'$ and $U$ will be the same, but the rim hooks will end in adjacent columns, so $\sigma_{\varrho(U')}+\sigma_{\varrho(U)}=0$. And we will have $(U')'=U$. It will follow that $\sum\sigma_T=0$, the sum over all $T$ descended from $U$ or $U'$. The proof of the proposition concludes by an application of Lemma 3 of \S 3.

There are two cases, $(a)$ and $(b)$, which are interchanged by the involution
$U\longleftrightarrow U'$. Let $V=U_{i-1}$, with shape $\alpha/\lambda$ (with $\alpha=\lambda$ if $i=1$). Case $(a)$ occurs when the $n$-rim in $\varrho(U)$
ends in a box of $\alpha$, and $(b)$ occurs when it ends in a box of $\varrho(U)/\alpha$. If $U$ is as in case $(a)$, the $n$-rim ends in a column $t$
with $t<k$, since $U$ is not proper; $U'$ is obtained from $U$ by taking 
$\widetilde{\varrho(U)}_t-\widetilde{\varrho(U)}_{t+1}$ $i$'s from the end of the first column and putting them at the end of column $t+1$. (That this many $i$'s can be removed from the first column follows from the fact that $V$ is proper.) Conversely, if $U$ is in case $(b)$, and the $n$-rim ends in column $t+1$, the $i$'s in the $n$-rim in column $t+1$ are removed and placed in the first column.\qed

\vskip 12pt

It follows from Proposition 4.1 that for $\lambda,\mu,\nu\subset l\times k$, with $|\lambda|+|\mu|=|\nu|+mn,\; m\geq 0$,
\begin{equation}\label{qlrdual}
N_{\lambda\mu}^{\nu}(l,k)=
N_{\widetilde{\lambda}\widetilde{\mu}}^{\widetilde{\nu}}(k,l).
\end{equation}

Similarly, from Proposition 4.3, 
\begin{equation}\label{qkdual}
K_{\lambda\mu}^{\nu}(l,k)=
\widetilde{K}_{\widetilde{\lambda}\mu}^{\widetilde{\nu}}(k,l),
\end{equation}
for $\lambda,\nu\subset l\times k$, and
$\mu =(\mu_1,\ldots ,\mu_r),\; 0\leq\mu_i\leq k$.
The numbers on the right are the conjugate quantum Kostka numbers, calculated by regarding $\widetilde{\lambda}$ and $\widetilde{\nu}$ in the $k\times l$ rectangle.
Equivalently, if $\nu[\widetilde{m}]$ is the partition obtained from $\nu$ by adding $m$ rim $n$-hooks, each occupying rows $1$ through $l$, then 
$\widetilde{K}_{\widetilde{\lambda}\mu}^{\widetilde{\nu}}(k,l)$ is the number of tableaux $T$ of content $\mu$ on $\nu[\widetilde{m}]/\lambda$ such that for each
entry $i$ of $T$ in the first row and column $k+p,\; p\geq 1$, the box in the $l^{\rm th}$ row and column $p$ is either in $\lambda$, or it has an entry of $T$ strictly smaller than $i$. We do not know a correspondance between the sets of tableaux that would explain equation $(\ref{qkdual})$.

\begin{example}
{\rm $k=5$, $l=4$, $\lambda=(5,3,3,1)$, $\mu=(5,2,2)$, $\nu=(2,1)$:}

\vskip 10pt
\epsfxsize 9cm
\centerline{\epsfbox{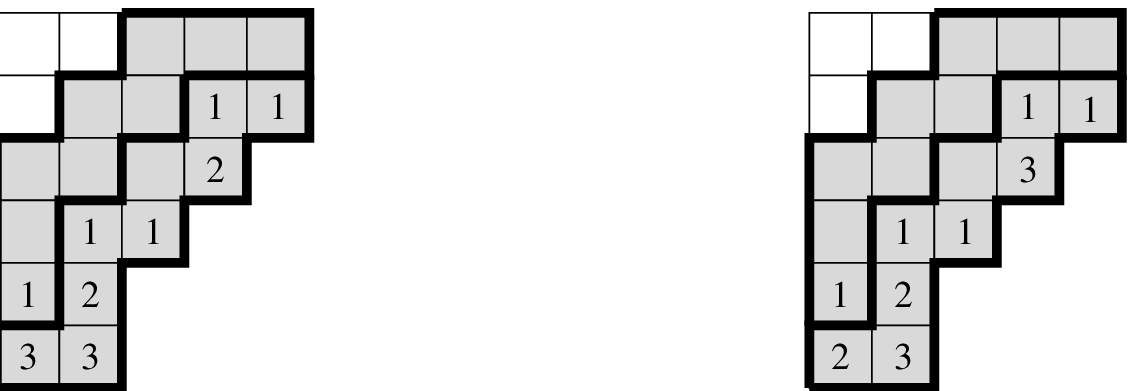}}
\vskip 10pt

\centerline{$K_{\lambda\mu}^{\nu}(4,5)=2$.}
\vskip 10pt
{\rm Note that there are in fact no nonproper tableaux of content $\mu$ on
the shape $\nu[m]/\lambda$.}

\vskip 10pt
\epsfxsize 8cm
\centerline{\epsfbox{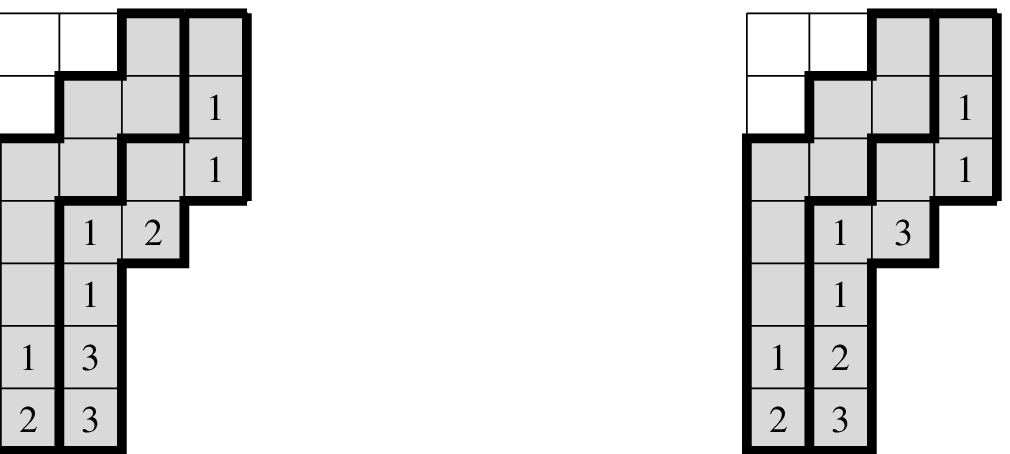}}
\vskip 10pt

\centerline{$\widetilde{K}_{\widetilde{\lambda}\mu}^{\widetilde{\nu}}(5,4)=2$.}

\vskip 10pt
{\rm As opposed to the previous case, there are many nonproper conjugate tableaux on $\widetilde{\nu}[m]/\lambda$. One such tableau is shown below.}

\vskip 10pt
\epsfxsize 2cm
\centerline{\epsfbox{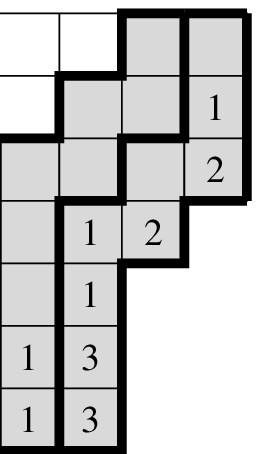}}
\vskip 10pt
\end{example}

Since
$$K_{\lambda\mu}^{\nu}(l,k)=\sum\epsilon(\varrho/\nu)K_{\lambda\mu}^{\varrho},$$
where the sum is over all $\varrho$ that can be obtained from $\nu$ by adding $m$ rim $n$-hooks, and
$K_{\lambda\mu}^{\varrho}$ is the classical Kostka number, equation $(\ref{qkdual})$ gives a nontrivial relation among classical Kostka numbers.

\end{document}